\documentclass[journal,final]{IEEEtran}

\usepackage{cite}
\usepackage{graphicx}
\usepackage{epsfig}
\usepackage{epstopdf}
\usepackage{amsmath}
\usepackage{amsfonts}
\usepackage{psfrag}
\usepackage{rotating}
\usepackage{latexsym}
\usepackage{dsfont}
\usepackage{stmaryrd}
\usepackage{amsthm}
\usepackage{amssymb}
\usepackage{fancyvrb}
\usepackage{textcomp}
\usepackage{url}
\usepackage{ifdraft}
\usepackage{url}
\usepackage{multirow}
\usepackage{hyperref}
\usepackage{rotating}
\usepackage{xspace}
\usepackage{array}
\usepackage[usenames]{color}
\usepackage{arydshln}
\usepackage{multido}
\usepackage{flushend}
\usepackage{enumerate}
\usepackage{caption}
\usepackage{subcaption}
\usepackage[latin1]{inputenc}
\usepackage{makecell}
\graphicspath{{figure/}{/}}

\newtheorem{specialcasecounter}{Theorem}
\newtheoremstyle{specialcasestyle}{1mm}{1mm}{\upshape}{}{\bfseries\upshape}{.}{0mm}{}
\theoremstyle{specialcasestyle}

\newcommand{\emptyline}{\raisebox{0.8ex}{\rule{24pt}{0.4pt}}}

\newcommand{\figref}[1]{Fig.~\protect\ref{#1}}

\newtheorem{corollary}{Corollary}[specialcasecounter]

\begin{document}

\title{New Results On the Sum of Two Generalized Gaussian Random Variables}

\author{Hamza~Soury,~\IEEEmembership{Student Member,~IEEE}, and~Mohamed-Slim~Alouini,~\IEEEmembership{Fellow,~IEEE}
\thanks{H.~Soury and M.-S.~Alouini are with King Abdullah University of Science and Technology (KAUST), Thuwal, Makkah Province, Kingdom of Saudi Arabia. The authors are members of the KAUST Strategic Research Initiative (SRI) in Uncertainty Quantification. E-mail:~\{soury.hamza, slim.alouini\}@kaust.edu.sa.}
	\vspace{-5mm}	
}

\maketitle

\begin{abstract}
We propose in this paper a new method to compute the characteristic function (CF) of generalized Gaussian (GG) random variable in terms of the Fox H function. The CF of the sum of two independent GG random variables is then deduced. Based on this results, the probability density function (PDF) and the cumulative distribution function (CDF) of the sum distribution are obtained. These functions are expressed in terms of the bivariate Fox H function. Next, the statistics of the distribution of the sum, such as the moments, the cumulant, and the kurtosis, are analyzed and computed. Due to the complexity of bivariate Fox H function, a solution to reduce such complexity is to approximate the sum of two independent GG random variables by one GG random variable with suitable shape factor. The approximation method depends on the utility of the system so three methods of estimate the shape factor are studied and presented.
\end{abstract}


\IEEEpeerreviewmaketitle

\section{Introduction}
The generalized Gaussian (GG) signals catched the interest of researches in the recent years. The GG distribution (GGD) is widely used to  model additive noise and interference, since several perturbations may indeed have a non-Gaussian behavior. It should be noted that the GGD family has a symmetric unimodal density characteristic with variable tail length that depends on the distribution parameter, the so called shape parameter $\alpha$.

In \cite{detectors_miller}, one of the non-Gaussian families used to model the noise is the GG noise (GGN). In fact, the GGN is deployed to study the detection for discrete time signals in non Gaussian noise. Moreover, in \cite{optimumdetection_yang,underwater_agrawal,novel_salahat} the noise is modeled by GGN and the performance analysis of the communication system perturbed by such noise has been studied and different metrics (such as probability of error, bit error rate, probability of false alarm etc.) have been evaluated. More recent works on the performance of communication systems perturbed by GGN can be found in \cite{average_soury,soury_error,exact_soury}, where the authors derived a closed form of the symbol/bit error rates of several modulations perturbed by GGN over a generalized flat fading channel. Furthermore, in \cite{simple_fiorina} the multi user interference in Ultra Wide Band was modeled as GGN. In addition, Fiorina introduced in\cite{simple_fiorina} a method to estimate the distribution parameter $\alpha$ using an estimation of the kurtosis of the GGD that depends only on the second and fourth moments.

From a mathematical point of view, let $X$ be a random variable following a GGD with parameter $\alpha>0$, mean $\mathbb{E}[X]=\mu$, and variance $\mathbb{E}[(X-\mu)^2]=\sigma^2$, i.e. $X\sim GGD(\mu,\sigma,\alpha)$. The PDF of $X$, noted $f_X(x)$, is defined in \cite[Eq. 2]{average_soury} as
\begin{equation}
\label{pdf_ggd}
f_X(x)=\frac{\alpha\Lambda}{2\Gamma(1/\alpha)}\exp\left(-\Lambda^{\alpha}|x-\mu|^{\alpha}\right)~\forall x\in\mathbb{R},
\end{equation}
where $\Lambda=\frac{\Lambda_0}{\sigma}=\frac{1}{\sigma}\sqrt{\frac{\Gamma(3/\alpha)}{\Gamma(1/\alpha)}}$ is the normalizing coefficient with respect to the shaping parameter $\alpha$. $\mathbb{E}$ denotes the expectation operator.

Note that the GGD can model different type of distribution, for example the Gaussian distribution is obtained by setting $\alpha$ equal to $2$, while for $\alpha=1$ one get the Laplacian distribution, and for $\alpha\to \infty $ we get the uniform distribution.

The complementary CDF (CCDF) of standard GGD (zero mean and unit variance) was given in \cite[Eq. A.4]{average_soury} and defined as the generalized Q-function $Q_{\alpha}(\cdot)$. It is given by
\begin{equation}
\label{ccdf_ggd}
Q_{\alpha}(x)=\frac{1}{2\Gamma(1/\alpha)}\Gamma(1/\alpha,\Lambda_0^{\alpha}x^{\alpha})~~\forall x\geq0,
\end{equation}
where $\Gamma(a,x)$ is the upper incomplete gamma function \cite[Eq. 6.5.3]{abramowitz}. Note here that $Q_{\alpha}(x)=1-Q_{\alpha}(-x)$ for negative argument. Once the CCDF of standard GGD defined, the CDF of X, noted as $F_X(x)$ can be expressed as
\begin{equation}
\label{cdf_ggd}
F_X(x)=1-Q_{\alpha}\left(\frac{x-\mu}{\sigma}\right).
\end{equation}
Other statistics of the GGD are investigated in what follows.

From the previous state of art, we notice that the majority of works focusing on the GGD. Actually, the works on the sum of two independent generalized Gaussian random variables (GGRV) are not as many as those related to a single GGRV. In fact, the most recent work on the properties of the sum of two independent GGRV is given in \cite{qian_sumgg}, where Zhao \textit{et al.} studied the properties of the PDF of the sum of two GGRV. They proved that such PDF has the same properties of the PDF of GGD (symmetry, convexity, monotonicity...) but they did not compute the PDF of the sum. On the other hand, they gave an approximation of the PDF of the sum of two identically independent distributed (i.i.d.) GGD. In this work, we are investigating the PDF of the sum of two independent not necessarily identically distributed GGRVs and focus on finding a closed form for such a PDF.

\subsection{Motivation and Applications}
As mentioned in \cite{qian_sumgg}, we need sometimes to study the sum of independent GGRV (SGG), or a linear system of GG white noise such as AR(1) process driven by a GG process. Furthermore, in seismic signal processing, the received seismic signal is modeled by a convolution of seismic signals, where seismic reflective coefficient can be modeled by GGD \cite{7zhao}. Thereby the distribution of the sum, that appears also as the convolution of the single distribution, is needed in this seismic model. Moreover, in communications, it is shown above that in some instances the noise and the multi user interference can be modeled as GG white noise \cite{detectors_miller,ultrawide_ahmed_park}. Therefore, the total perturbations at reception, defined as the sum of noise and interference, is modeled as the sum the GG signals. Many other applications can be found in the literature to motivate the study of sum distribution of GG signals.

Per consequence, it is important to study the statistics and the density of the SGG distribution. The PDF was not derived before and so the CDF. Actually, the approach used in this work is based on the CF which was investigated in \cite{tibor_function} for $\alpha>1$. However in our case we are studying the CF for any value of $\alpha$ using another approach of calculation based on the properties of the Fox H function (FHF) \cite{h_function}. From the CF of one GGRV, we get the CF of the sum. Hence, the relation between the CF and the distribution densities leads us to investigate the PDF and the CDF of such distribution and its statistics (moments, cumulant, kurtosis...).

\subsection{Organization of the Paper}
The remaining of this paper is organized as follows. In Section II, we investigate the statistics of the GGD, such as the CF, the moments, cumulant, and kurtosis. In Section III, the PDF and some statistics of the SGG distribution are presented. Next in Section IV, we analyze a method to approximate the PDF of the sum by a PDF of one GGD and the performance of such approximation is studied. Finally Section V concludes the paper.

\section{Generalized Gaussian Distribution Statistics}
\subsection{Characteristic Function}
Let $\alpha>0$, and $X$ be random variable (RV) following a $GGD(\mu,\sigma,\alpha)$
\begin{specialcasecounter}
The CF of $X$, $\mathbb{E}[e^{itX}]$, is given by
\begin{align}
\varphi_{\alpha}(t)=\frac{\sqrt{\pi}}{\Gamma(1/\alpha)}&e^{it\mu}
		{\rm{H}}_{1,2}^{1,1}
			\left[
				\frac{\sigma^2\Gamma(1/\alpha)}{4\Gamma(3/\alpha)}t^2
				\left|
					\genfrac{}{}{0pt}{}
						{(1-\frac{1}{\alpha},\frac{2}{\alpha})}
						{(0,1),(\frac{1}{2},1)}
				\right.
			\right],
\end{align}
\end{specialcasecounter}
where $\rm{H}^{\cdot,\cdot}_{\cdot,\cdot}[\cdot]$ is the Fox H function (FHF) \cite[Eq. (1.1.1)]{h_function},\cite{h_function_int}.
\begin{proof}
Let us start from the definition of the CF and the PDF of the GGD. In fact the CF can be expressed in integral form as
\begin{align}
\label{cf_integral1}
\varphi_{\alpha}(t)&=&\mathbb{E}[e^{itX}]=&\int_{\mathbb{R}}e^{itx}f_X(x)dx\nonumber\\
&=&\frac{\alpha\Lambda}{2\Gamma(1/\alpha)}&e^{it\mu}\int_{\mathbb{R}}e^{itx}\exp\left(-\Lambda^{\alpha}|x|^{\alpha}\right)dx.
\end{align}
Since $|x|$ is an even function, the integral in \eqref{cf_integral1} is the cosine transform of the exponential component
\begin{align}
\label{cf_integral2}
\varphi_{\alpha}(t)=\frac{\alpha\Lambda}{\Gamma(1/\alpha)}e^{it\mu}\int_{0}^{\infty}\cos(tx)\exp\left(-(\Lambda x)^{\alpha}\right)dx.
\end{align}
Alternative expressions of $\cos(x)$ and $e^{x^{\alpha}}$ in terms of the FHF are available in \cite[Eq. (2.9.8) \& Eq. (2.9.4)]{h_function}
\begin{align}
\label{cosine_H}
\cos(x)&=\sqrt{\pi}~{\rm{H}}_{0,2}^{1,0}\left[\frac{x^2}{4}\left|\genfrac{}{}{0pt}{}{\emptyline}
{(0,1),(\frac{1}{2},1)}\right.\right] \\
e^{-x^{\alpha}}&=\frac{1}{\alpha}{\rm{H}}_{0,1}^{1,0}\left[x\left|\genfrac{}{}{0pt}{}{\emptyline}
{(0,\frac{1}{\alpha})}\right.\right]
\end{align}
Hence the integral identity defined in \cite[Eq. (2.8.4)]{h_function}, solves the integral of the product of two FHF function over the positive real numbers. As a consequence the CF can be re-written as
\begin{align}
\label{proof_cf1}
\varphi_{\alpha}(t)&=\frac{\sqrt{\pi}\Lambda}{\Gamma(1/\alpha)}e^{it\mu}\int_{0}^{\infty}
{\rm{H}}_{0,2}^{1,0}\left[\frac{t^2}{4}x^2\left|\genfrac{}{}{0pt}{}{\emptyline}
{(0,1),(\frac{1}{2},1)}\right.\right]\times\nonumber\\
&\qquad\qquad\qquad\qquad{\rm{H}}_{0,1}^{1,0}\left[\Lambda x\left|\genfrac{}{}{0pt}{}{\emptyline}
{(0,\frac{1}{\alpha})}\right.\right] dx\nonumber \\
&=
\frac{\sqrt{\pi}}{\Gamma(1/\alpha)}e^{it\mu}
		{\rm{H}}_{1,2}^{1,1}
			\left[
				\frac{\sigma^2\Gamma(1/\alpha)}{4\Gamma(3/\alpha)}t^2
				\left|
					\genfrac{}{}{0pt}{}
						{(1-\frac{1}{\alpha},\frac{2}{\alpha})}
						{(0,1),(\frac{1}{2},1)}
				\right.
			\right].
\end{align}
\end{proof}
Note that the result proved in \eqref{proof_cf1} is valid for all positive shape parameter $\alpha>0$. A previous demonstration was derived in \cite{tibor_function} but only for $\alpha>1$, while the authors provided an expression of the CF of GGD in terms of the Fox-Wright generalized hypergeometric function $p\Psi q(\cdot)$, which is a special case of the FHF \cite[Eq. (2.9.29)]{h_function}. The CF expression and derivation presented in this work present one of the contributions of the paper.

\subsection{Moment Generating Function}
The moment generating function (MGF) can be directly concluded from the CF by the relation $M_{\alpha}(t)=\varphi_{\alpha}(-it)$, so the MGF is obtained by
\begin{equation}
\label{mgf_ggd}
M_{\alpha}(t)=\frac{\sqrt{\pi}e^{t\mu}}{\Gamma(1/\alpha)}
		{\rm{H}}_{1,2}^{1,1}
			\left[
				\frac{-\sigma^2\Gamma(1/\alpha)}{4\Gamma(3/\alpha)}t^2
				\left|
					\genfrac{}{}{0pt}{}
						{(1-\frac{1}{\alpha},\frac{2}{\alpha})}
						{(0,1),(\frac{1}{2},1)}
				\right.
			\right].
\end{equation}
In some special cases, the MGF can be expressed in terms of elementary functions. For example for Gaussian case, i.e. $\alpha=2$, and using the special case of the FHF in \cite[Eq. (2.9.4)]{h_function}, the MGF of Gaussian is
\begin{equation}
\label{mgf_gaussian}
M_{2}(t)=e^{t\mu}
		{\rm{H}}_{0,1}^{1,0}
			\left[
				\frac{-\sigma^2}{2}t^2
				\left|
					\genfrac{}{}{0pt}{}
						{\emptyline}
						{(0,1)}
				\right.
			\right]
=e^{t\mu+\frac{1}{2}\sigma^2t^2}
\end{equation}
Another special case appears interesting is the Laplacian distribution (i.e. $\alpha=1$). Using the identity \cite[Eq. (2.9.5)]{h_function}, the MGF of Laplacian distribution is given by
\begin{align}
\label{mgf_laplace}
M_{1}(t)&=\sqrt{\pi}e^{t\mu}
		{\rm{H}}_{1,2}^{1,1}
			\left[
				\frac{-\sigma^2}{8}t^2
				\left|
					\genfrac{}{}{0pt}{}
						{(0,2)}
						{(0,1),(\frac{1}{2},1)}
				\right.
			\right]\nonumber\\
&=e^{t\mu}
		{\rm{H}}_{1,1}^{1,1}
			\left[
				-\frac{1}{2}\sigma^2t^2
				\left|
					\genfrac{}{}{0pt}{}
						{(0,1)}
						{(0,1)}
				\right.
			\right]\nonumber \\
&=\frac{e^{t\mu}}{1-\frac{1}{2}\sigma^2t^2}.
\end{align}

\subsection{Moments and Cumulant}
Without loss of generality, we are focusing our study to zero mean random variables (i.e. $\mu=0$).

Due to the symmetry of the PDF of GGD, the odd moments vanish and the even moments obtained as follows
\begin{equation}
\label{moment_ggd}
\left\{
\begin{array}{lll}
m_{2n}(X)&=\mathbb{E}[X^{2n}]&=\sigma^{2n}\frac{\Gamma(1/\alpha)^n}{\Gamma(3/\alpha)^n}\frac{\Gamma(\frac{2n+1}{\alpha})}{\Gamma(1/\alpha)}\\
m_{2n+1}(X)&=\mathbb{E}[X^{2n+1}]&=0
\end{array}
\right.
\end{equation}
Once the MGF and the moments are obtained, one can investigate the expression of the cumulant generating function (CGF) and the cumulant of the GGD. Actually the CGF, $K_X(t)$ (or $K_{\alpha}(t)$), is defined as
\begin{align}
\label{cgf_ggd1}
K_{\alpha}(t)&=\log M_{\alpha}(t)\nonumber\\
&=\log(\frac{\sqrt{\pi}}{\Gamma(1/\alpha)})\nonumber\\
&\qquad+\log{\rm{H}}_{1,2}^{1,1}
			\left[
				\frac{-\sigma^2\Gamma(1/\alpha)}{4\Gamma(3/\alpha)}t^2
				\left|
					\genfrac{}{}{0pt}{}
						{(1-\frac{1}{\alpha},\frac{2}{\alpha})}
						{(0,1),(\frac{1}{2},1)}
				\right.
			\right].
\end{align}
By definition, the n-\textbf{th cumulant}, noted $k_{n}(X)$, is the n-th term in the Taylor series expansion of $K_{\alpha}(t)$ at $0$.
\begin{specialcasecounter}
The even cumulants of a zero mean GG random variable $X$ can be expressed in terms of the even moments of $X$ by
\begin{align}
\label{cumulant_ggd}
k_{2n}(X)&=-\sum_{m_1+2m_2+\dots n m_n=n}\frac{(2n)!(m_1+\dots+m_n-1)!}{m_1!m_2!\dots m_n!}\times\nonumber\\
&\prod_{1\leq j\leq n}\left(-\frac{\sigma^{2j}\Gamma(1/\alpha)^j\Gamma(\frac{2j+1}{\alpha})}{\Gamma(3/\alpha)^j\Gamma(1/\alpha)(2j)!}\right)^{m_j},
\end{align}
and the odd cumulants are equal to zero.
\end{specialcasecounter}
\begin{IEEEproof}
A cumulant $k_{n}(X)$ is the n-th derivative of the CGF evaluated at zero, $k_{n}(X)=\left.\frac{d^nK_{\alpha}(t)}{dt^n}\right|_{t=0}=\frac{d^nK_{\alpha}(0)}{dt^n}$. Since the CGF appears as the composite of two functions, we may use the Fa\`a di Bruno's formula \cite[Eq. (2)]{bruno_craik} that computes the n-th derivative of composite functions
\begin{align}
\label{nthderivative_cumulant}
\frac{d^nK_{\alpha}(t)}{dt^n}=-\sum_{m_1,\dots,m_n}\frac{n!(m_1+\dots+m_n-1)!}{m_1!m_2!\dots m_n!}\times\nonumber\\
\prod_{1\leq j\leq n}\left(-\frac{M_{\alpha}^{(j)}(t)}{j!M_{\alpha}(t)}\right)^{m_j},
\end{align}
the sum is over $m_1,m_2,\dots,m_n$ such that $\displaystyle\sum_{j=1}^n jm_j=n$. Thereby evaluating \eqref{nthderivative_cumulant} at zero and replacing the moment by its expression, one get the final expression of the cumulant.
\end{IEEEproof}
The cumulant of low order are easy to expressed
\begin{align}
\label{cumulant_loworder}
k_2(X)&=\sigma^2\nonumber\\
k_4(X)&=\sigma^4\left(\frac{\Gamma(1/\alpha)\Gamma(5/\alpha)}{\Gamma(3/\alpha)^2}-3\right)\nonumber\\
k_6(X)&=\sigma^6\left(\frac{\Gamma(1/\alpha)^2\Gamma(7/\alpha)}{\Gamma(3/\alpha)^3}-15\frac{\Gamma(1/\alpha)\Gamma(5/\alpha)}{\Gamma(3/\alpha)^2}+30\right).
\end{align}
Another statistics appear interesting to evaluate, in occurrence the kurtosis. The kurtosis, $Kurt(X)$, is defined as the fourth cumulant divided by the square of the second cumulant of the distribution. In the GGD case the kurtosis is equal to
\begin{align}
\label{kurtosis_ggd}
Kurt(X)=\frac{k_4(X)}{k_2(X)^2}=\frac{\Gamma(1/\alpha)\Gamma(5/\alpha)}{\Gamma(3/\alpha)^2}-3.
\end{align}
One can easily check that \eqref{kurtosis_ggd} confirms that the Gaussian kurtosis is equal to $0$ and the Laplacian kurtosis is equal to $3$.

At this stage, the statistics of one GGRV are expressed in closed form. The next section considers the densities and statistics of the SGG distribution.

\section{Sum of Two Independent GG Random Variables}
As known the CF of the sum of two independent RV is the product of their CFs. Since the CF of the GGD is defined in the previous section, the CF of the sum can be easily obtained and so the densities by inverse Laplace transform of the CF. In fact, let $X\sim GGD(\mu_1,\sigma_1,\alpha)$ and $Y\sim GGD(\mu_2,\sigma_2,\beta)$ two independent random variables following a GGD, and let $Z=X+Y$ be their sum. It is clear that the first and second moment of $Z$ are easy to find
\begin{align}
\label{mean_variance_sumz}
\mathbb{E}[Z]&=\mu=\mu_1+\mu_2\nonumber\\
\mathbb{E}[(Z-\mu)^2]&=\sigma^2=\sigma_1^2+\sigma_2^2.
\end{align}

\subsection{PDF and CDF of the Sum of Two GGRV}
The PDF of a random variable is known as the inverse Laplace transform of the CF. The CF of $Z$ is given by
\begin{align}
\label{cf_sumz}
\varphi_Z(t)&=\varphi_X(t)\varphi_Y(t)\nonumber\\
&=\frac{\pi e^{it\mu}}{\Gamma(1/\alpha)\Gamma(1/\beta)}
		{\rm{H}}_{1,2}^{1,1}
			\left[
				\frac{\sigma_1^2\Gamma(1/\alpha)}{4\Gamma(3/\alpha)}t^2
				\left|
					\genfrac{}{}{0pt}{}
						{(1-\frac{1}{\alpha},\frac{2}{\alpha})}
						{(0,1),(\frac{1}{2},1)}\right.\right]\times\nonumber\\
&\qquad\qquad\qquad		{\rm{H}}_{1,2}^{1,1}\left[
				\frac{\sigma_2^2\Gamma(1/\beta)}{4\Gamma(3/\beta)}t^2\left|
					\genfrac{}{}{0pt}{}
						{(1-\frac{1}{\beta},\frac{2}{\beta})}
						{(0,1),(\frac{1}{2},1)}\right.\right].
\end{align}
By applying the Laplace transform inverse to \eqref{cf_sumz}, the PDF of $Z$ is given by the following theorem.
\begin{figure*}
\vspace{-2mm}
\begin{minipage}[t]{0.9775\textwidth}
\begin{center}
\begin{align}
\label{pdf_sumz}
f_Z(z)&=\frac{\sqrt{\pi}}{\Gamma(1/\alpha)\Gamma(1/\beta)|z-\mu|} \times\nonumber\\
&\qquad\qquad{\rm{H}}^{0,1;1,1;1,1}_{2,0;1,2;1,2}\left[
\frac{\sigma_1^2\Gamma(1/\alpha)}{\Gamma(3/\alpha)(z-\mu)^2},\frac{\sigma_2^2\Gamma(1/\beta)}{\Gamma(3/\beta)(z-\mu)^2}\left|
\genfrac{}{}{0pt}{}{(\frac{1}{2},1,1),(0,1,1)}{\emptyline}\right| \left. \genfrac{}{}{0pt}{}{(1-\frac{1}{\alpha},\frac{2}{\alpha})}{(0,1),(\frac{1}{2},1)}
\right|\genfrac{}{}{0pt}{}{(1-\frac{1}{\beta},\frac{2}{\beta})}{(0,1),(\frac{1}{2},1)}\left.
\right.
\right]
\end{align}
\begin{align}
\label{cdf_sumz}
F_Z(z)&=\frac{1}{2}+\frac{\sqrt{\pi}~sign(z-\mu)}{2\Gamma(1/\alpha)\Gamma(1/\beta)} \times\nonumber\\
&\qquad\qquad{\rm{H}}^{0,1;1,1;1,1}_{2,0;1,2;1,2}\left[
\frac{\sigma_1^2\Gamma(1/\alpha)}{\Gamma(3/\alpha)(z-\mu)^2},\frac{\sigma_2^2\Gamma(1/\beta)}{\Gamma(3/\beta)(z-\mu)^2}\left|
\genfrac{}{}{0pt}{}{(\frac{1}{2},1,1),(1,1,1)}{\emptyline}\right| \left. \genfrac{}{}{0pt}{}{(1-\frac{1}{\alpha},\frac{2}{\alpha})}{(0,1),(\frac{1}{2},1)}
\right|\genfrac{}{}{0pt}{}{(1-\frac{1}{\beta},\frac{2}{\beta})}{(0,1),(\frac{1}{2},1)}\left.
\right.
\right]
\end{align}
\end{center}
\vspace{1pt}
\end{minipage}\\
\rule{\textwidth}{0.4pt}
\noindent\vspace{-15pt}
\end{figure*}

\begin{specialcasecounter}
\label{pdf_theorem}
The PDF of the sum of two independent GG random variable can be expressed in terms of the bivariate FHF \cite{mittal_biv} in \eqref{pdf_sumz}.
\end{specialcasecounter}

The FHF of two variables \cite{mittal_biv}, also known as the Bivariate Fox H-function (BFHF) $\rm{H}_{\cdot~\cdot;\cdot~\cdot;\cdot~\cdot}^{\cdot~\cdot;\cdot~\cdot;\cdot~\cdot}[\cdot,\cdot]$ is a generalization of the FHF. Its MATLAB implementation is outlined in\cite{peppas_formula}.
\begin{IEEEproof}
Let $A=\frac{\sigma_1^2\Gamma(1/\alpha)}{4\Gamma(3/\alpha)}$ and $B=\frac{\sigma_2^2\Gamma(1/\beta)}{4\Gamma(3/\beta)}$. The inverse Laplace transform of the CF \eqref{cf_sumz} of $Z$ gives the PDF
\begin{align}
\label{pdf_integral1}
f_Z(z)&=\frac{1}{2\pi}\int_{\mathbb{R}}e^{-itz}\varphi_Z(t)dt\nonumber\\
&=\frac{1}{2\Gamma(1/\alpha)\Gamma(1/\beta)}
\int_{\mathbb{R}}
		{\rm{H}}_{1,2}^{1,1}
			\left[
				A~t^2
				\left|
					\genfrac{}{}{0pt}{}
						{(1-\frac{1}{\alpha},\frac{2}{\alpha})}
						{(0,1),(\frac{1}{2},1)}\right.\right]\times\nonumber\\
&\qquad\qquad		{\rm{H}}_{1,2}^{1,1}\left[
				B~t^2\left|
					\genfrac{}{}{0pt}{}
						{(1-\frac{1}{\beta},\frac{2}{\beta})}
						{(0,1),(\frac{1}{2},1)}\right.\right]e^{it(\mu-z)}dt
\end{align}
The first two FHF are even functions, so the integral becomes a cosine transform of the product of these two FHF functions
\begin{align}
\label{pdf_integral2}
f_Z(z)&=\frac{1}{2\Gamma(1/\alpha)\Gamma(1/\beta)}
\int_{0}^{\infty}
		{\rm{H}}_{1,2}^{1,1}
			\left[
				A~t^2
				\left|
					\genfrac{}{}{0pt}{}
						{(1-\frac{1}{\alpha},\frac{2}{\alpha})}
						{(0,1),(\frac{1}{2},1)}\right.\right]\times\nonumber\\
&\qquad\qquad		{\rm{H}}_{1,2}^{1,1}\left[
				B~t^2\left|
					\genfrac{}{}{0pt}{}
						{(1-\frac{1}{\beta},\frac{2}{\beta})}
						{(0,1),(\frac{1}{2},1)}\right.\right]\cos\left(t(\mu-z)\right)dt
\end{align}
As seen before, the cosine has a representation in terms of the FHF \eqref{cosine_H}. So we are facing an integral that involves the product of three FHFs over the positive real numbers. Such integral is solved in\cite[Eq.~(2.3)]{mittal_biv} and it is expressed in terms of the BFHF, which give us the final expression of the PDF of $Z$.
\end{IEEEproof}
The CDF of $Z$ is the primitive of $f_Z$ that vanishes at $(-\infty)$. Back to \eqref{pdf_integral2}, it appears that the CDF is expressed in term of an integral involving the product of two FHFs and sine function. The latter can be expressed in terms of the FHF for positive argument \cite[Eq. (2.9.7)]{h_function}. Thereby, the CDF of the sum of two independent GGRV becomes the integral of the product of three FHFs which is evaluated in terms of the BFHF.
\begin{corollary}
\label{cdf_theorem}
The CDF, $F_Z(z)$, of the SGG distribution is given in \eqref{cdf_sumz}.
\end{corollary}
In \eqref{cdf_sumz}, $sign(x)$ gives the sign of the real number $x$. The results in Theorem.\ref{pdf_theorem} and Corollary \ref{cdf_theorem} represent new results and they were not investigated before, which make them the essential contribution of this paper.

\subsection{Statistics of $Z$}
In the following analysis, the zero mean case is considered, while the non zero mean random variable can be obtained from the zero mean random variable by a simple shift $Z=Z_0+\mu$.

\subsubsection{MGF}
As mentioned before the MGF of $Z$ can be obtained from the CF by the relation $M_Z(t)=\varphi_Z(-it)$ which gives the MGF of Z as
\begin{align}
\label{mgf_sumz}
M_Z(t)&=\frac{\pi}{\Gamma(1/\alpha)\Gamma(1/\beta)}
		{\rm{H}}_{1,2}^{1,1}
			\left[
				\frac{-\sigma_1^2\Gamma(1/\alpha)}{4\Gamma(3/\alpha)}t^2
				\left|
					\genfrac{}{}{0pt}{}
						{(1-\frac{1}{\alpha},\frac{2}{\alpha})}
						{(0,1),(\frac{1}{2},1)}\right.\right]\nonumber\\
&\qquad\qquad\qquad		{\rm{H}}_{1,2}^{1,1}\left[
				\frac{-\sigma_2^2\Gamma(1/\beta)}{4\Gamma(3/\beta)}t^2\left|
					\genfrac{}{}{0pt}{}
						{(1-\frac{1}{\beta},\frac{2}{\beta})}
						{(0,1),(\frac{1}{2},1)}\right.\right].
\end{align}

\subsubsection{Moment} The moments of $Z$ can be obtained from the binomial formula that describes the integer power of the sum of two numbers. Hence known the moment of $X$ and $Y$ \eqref{moment_ggd}, it appears that the odd moments of $Z$, $m_{2n+1}(Z)$, vanish while the even moments are given by
\begin{align}
\label{moment_sumz}
m_{2n}(Z)&=\frac{\sigma_2^{2n}\Gamma(\frac{1}{\beta})^n}{\Gamma(\frac{1}{\alpha})\Gamma(\frac{1}{\beta})\Gamma(\frac{3}{\beta})^n}\sum_{k=0}^{n} {2n \choose 2k}\left(\frac{\sigma_1^{2}}{\sigma_2^{2}}\frac{\Gamma(\frac{1}{\alpha})\Gamma(\frac{3}{\beta})}{\Gamma(\frac{3}{\alpha})\Gamma(\frac{1}{\beta})}\right)^k\times\nonumber\\
&\qquad\qquad\qquad\Gamma\left(\frac{2k+1}{\alpha}\right)\Gamma\left(\frac{2n-2k+1}{\beta}\right)
\end{align}

\subsubsection{Cumulant and Kurtosis}
From the MGF, it is easier to get the CGF by applying the logarithm to the MGF $M_Z(t)$. Thereby the CGF of $Z$ is the sum of the CGF of $X$ and the CGF of $Y$, $K_Z(t)=K_X(t)+K_Y(t)$. Moreover, the cumulant of $Z$ is expressed also as the sum of the cumulant of $X$ and the cumulant of $Y$. Note that the odd cumulant are equal to zero while the even ones are given by
\begin{align}
\label{cumulant_sumz}
k_{2n}(Z)&=k_{2n}(X)+k_{2n}(Y)\nonumber\\
&=-\!\!\!\!\sum_{m_1+2m_2+\dots n m_n=n}\!\!\!\!\frac{(2n)!(m_1+\dots+m_n-1)!}{m_1!m_2!\dots m_n!}\times\nonumber\\
&\left[\prod_{1\leq j\leq n}\left(-\frac{\sigma_1^{2j}\Gamma(\frac{1}{\alpha})^j\Gamma(\frac{2j+1}{\alpha})}{\Gamma(\frac{3}{\alpha})^j\Gamma(\frac{1}{\alpha})(2j)!}\right)^{m_j}\right.\nonumber\\
&\left.+
\prod_{1\leq k\leq n}\left(-\frac{\sigma_2^{2k}\Gamma(\frac{1}{\beta})^k\Gamma(\frac{2k+1}{\beta})}{\Gamma(\frac{3}{\beta})^k\Gamma(\frac{1}{\beta})(2k)!}\right)^{m_k}
\right].
\end{align}
Once the cumulant expression is evaluated, the kurtosis of $Z$ can be expressed, per definition, in terms of the fourth and second moments as
\begin{align*}
Kurt(Z)=\frac{k_4(Z)}{k_2(Z)^2}=\frac{k_4(X)+k_4(Y)}{\left(k_2(X)+k_2(Y)\right)^2}.
\end{align*}
Note that $k_2(X)+k_2(Y)=\sigma_1^2+\sigma_2^2=\sigma^2$. Thus a relation between the kurtosis of $Z$, $X$ and $Y$ appears as
\begin{align}
\label{kurtosis_sumz}
Kurt(Z)=\frac{\sigma_1^2}{\sigma^2}Kurt(X)+\frac{\sigma_2^2}{\sigma^2}Kurt(Y).
\end{align}
The final expression of the kurtosis of $Z$ is thus given by
\begin{align}
\label{kurtosis_sumz}
&Kurt(Z)=\nonumber\\
&\qquad\frac{\sigma_1^4}{\sigma^4}\frac{\Gamma(\frac{1}{\alpha})\Gamma(\frac{5}{\alpha})}{\Gamma(\frac{3}{\alpha})^2}
+\frac{\sigma_2^4}{\sigma^4}\frac{\Gamma(\frac{1}{\beta})\Gamma(\frac{5}{\beta})}{\Gamma(\frac{3}{\beta})^2}+6\frac{\sigma_1^2\sigma_2^2}{\sigma^4}-3
\end{align}

\section{Approximation of the PDF of the sum of two GGRV}
The expression of the PDF of the sum of two independent GGRV \eqref{pdf_sumz} is quite high complex since it is expressed in terms of the BFHF. Therefore, an approximation of the PDF is highly recommended to simplify the calculations and study, in simple way, the performance of systems in which the PDF of the sum is needed, like, for example, the evaluation of the symbol error rate (SER) of an M-phase shift keying (MPSK) over an GGN channel. Such analysis needs the PDF and the CDF of the SGG distribution.

In this section we are investigating the approximation of the PDF of $Z$ by the PDF of another GG random variable with shape factor $\gamma$ to be determined. In \cite{qian_sumgg}, it has been proved that the PDF of the sum cannot be a PDF of one GGRV. However the authors proved that both PDFs have the same properties (symmetric, convexity, monotonicity...). Furthermore, the PDF of the sum of two i.i.d. GGRV was approximated by the PDF of GGD. From that analysis, an approximation of the PDF of $Z$ by the PDF of GGD is needed and worth pursuing.

As shown in \eqref{pdf_ggd}, $3$ parameters are needed to characterize a GGD, namely, the mean, the variance and the shape factor. The mean and the variance are given in \eqref{mean_variance_sumz}. Therefore, we need to find a method to get the shape factor $\gamma$. In what follows, three approaches are presented.

\subsection{Kurtosis Approach}
The first method to estimate $\gamma$ is by using the kurtosis of the distributions. Since the kurtosis of the sum is already known, the shape factor can be obtained by equalizing both kurtosis. Thereby, we get the following equation to solve
\begin{equation}
\label{gamma_kurtosis}
\begin{array}{c}
\displaystyle Kurt(Z_{\gamma})=Kurt(Z)\\
\displaystyle\iff\\
\displaystyle\frac{\Gamma(\frac{1}{\gamma})\Gamma(\frac{5}{\gamma})}{\Gamma(\frac{3}{\gamma})^2}=\frac{\sigma_1^4}{\sigma^4}\frac{\Gamma(\frac{1}{\alpha})\Gamma(\frac{5}{\alpha})}{\Gamma(\frac{3}{\alpha})^2}
+\frac{\sigma_2^4}{\sigma^4}\frac{\Gamma(\frac{1}{\beta})\Gamma(\frac{5}{\beta})}{\Gamma(\frac{3}{\beta})^2}+6\frac{\sigma_1^2\sigma_2^2}{\sigma^4},
\end{array}
\end{equation}
while $Z_{\gamma}\sim GGD(\mu,\sigma,\gamma)$ is the approximated RV of $Z$ with parameter $\gamma$.

\begin{figure}[h]
\begin{center}
\includegraphics[width=0.95\columnwidth]{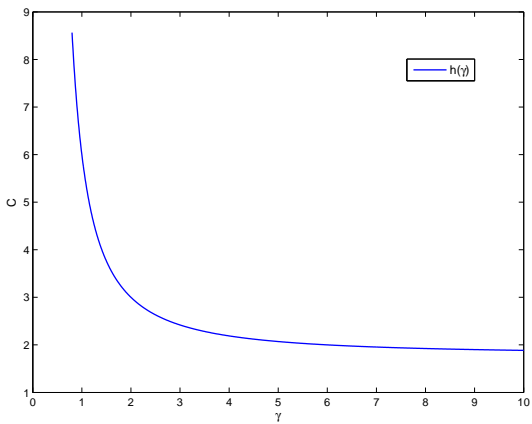}
\caption{The curve of $h(\gamma)$ for positive values of $\gamma$.}
\label{figure1}
\end{center}
\end{figure}

Lets define the ratio between the variance of $X$ and $Y$ as $\delta=\frac{\sigma_1^2}{\sigma_2^2}$, so the equation on $\gamma$ can be written in terms of $\delta$ as
\begin{equation}
\label{gamma_kurtosis1}
\begin{array}{c}
\displaystyle\frac{\Gamma(\frac{1}{\gamma})\Gamma(\frac{5}{\gamma})}{\Gamma(\frac{3}{\gamma})^2} =\frac{1}{(1+\delta)^2}\left(\delta^2\frac{\Gamma(\frac{1}{\alpha})\Gamma(\frac{5}{\alpha})}{\Gamma(\frac{3}{\alpha})^2}
+\frac{\Gamma(\frac{1}{\beta})\Gamma(\frac{5}{\beta})}{\Gamma(\frac{3}{\beta})^2}+6\delta\right),
\end{array}
\end{equation}
By knowing $\alpha$, $\beta$, and $\delta$, \eqref{gamma_kurtosis} is written as $h(\gamma)=C$, where $h(\gamma)$ is a function on $\gamma$, and $C$ is a known positive constant. In \figref{figure1}, the function $h$ is drawn versus $\gamma$ to analyze its behavior. Therefore, it appears that the function $h(\cdot)$ is a bijection. As such the equation $h(\gamma)=C$ has only one solution in the positive real axis. Which mean that $\gamma$ exists and it is unique. The value of $\gamma$ is given as $\gamma_{Kurt}$ in Table \ref{table1} for some scenarios along with other values of $\gamma$ obtained from other approaches that we will discuss later on.

\subsection{Best Tail Approximation}
Another method to estimate $\gamma$ consist of taking the best choice of $\gamma$ that minimizes the square error of the tail. In other words $\gamma$ is chosen so the error between the exact PDF, $f_Z(z)$ and the approximated PDF $f_{\gamma}(z)$ at the tail is minimal. The tail is defined so $z$ is above some level $z\geq n \sigma$
\begin{align}
\label{gamma_minimum_squareerror}
\gamma_{Tail}=\text{arg}\min_{\gamma>0}\int_{n\sigma}^{\infty}\left(f_{\gamma}(z)-f_Z(z)\right)^2dz,
\end{align}
where $n$ is chosen to define the desired region of the tail of the distribution. The minimization in \eqref{gamma_minimum_squareerror} cannot be solved analytically by the available tools since it contains a shifted integral of a BFHF and FHF which is not known yet. A numerical evaluation of $\gamma_{Tail}$ is given in Table \ref{table1} for different values of $\delta$.

\subsection{CDF Approximation}
This method is used to obtain the shape parameter that minimizes the error between the CDF of $Z$ and the approximated CDF. Such  approximation will give an asymptotic approximation of the complementary CDF (CCDF) which is needed in the computation of the probability of error. Mathematically, the shape parameter is given by
\begin{align}
\label{gamma_minimum_cdf}
\gamma_{CDF}=\text{arg}\min_{\gamma>0}\int_{0}^{\infty}\left(F_{\gamma}(z)-F_Z(z)\right)^2dz.
\end{align}
In Table \ref{table1} some numerical values of $\gamma_{CDF}$ are given and a comparison between three methods of shape parameter estimation is available too.

\begin{table}[!h]
\caption{Shape parameter for the approximated PDF using kurtosis, minimum CDF error, and minimum tail error for $\sigma_1=1$ and different values of $(\alpha,~\beta,~\delta)$}
\label{table1}
\begin{center}
\begin{tabular}{|c|c|c|c|c|c|c|}
\hline
\multirow{2}{*}{$(\alpha,\beta,\delta)$}&\multirow{2}{*}{$\gamma_{Kurt}$}&\multirow{2}{*}{$\gamma_{CDF}$}&
\multicolumn{4}{ c| }{$\gamma_{Tail}$}\\ \cline{4-7}
&&&$n=0$&$1$&$2$&$3$
\\
\hline
$(0.5,0.5,1)$&$0.626$&0.467&$0.768$&$0.673$&$0.624$&$0.642$\\
\hline
$(0.5,0.5,2)$&$0.604$&0.492&$0.762$&$0.656$&$0.603$&$0.584$\\
\hline
$(0.5,0.7,2)$&$0.633$&0.501&$0.861$&$0.741$&$0.636$&$0.834$\\
\hline
$(0.5,1.2,1)$&$0.779$&0.602&$1.160$&$1.053$&$0.757$&$1.165$\\
\hline
$(1.5,1.5,2)$&$1.673$&1.373&$1.738$&$1.702$&$1.683$&$1.664$\\
\hline
$(1.5,2.5,1)$&$1.908$&1.391&$1.979$&$1.959$&$1.952$&$1.887$\\
\hline
$(1.5,2.5,2)$&$1.753$&1.443&$1.842$&$1.799$&$1.771$&$1.741$\\
\hline
$(2.5,~3~,3)$&$2.295$&1.941&$2.226$&$2.261$&$2.267$&$2.335$\\
\hline
\end{tabular}
\end{center}
\end{table}
An overview from Table \ref{table1} shows that the optimal value of $\gamma_{Tail}$ is near the value given by the kurtosis for any value of $n$. It is clear also that $\gamma_{Tail}$ approaches closely to $\gamma_{Kurt}$ specially for $n=2$ for all values of $\delta$. This analysis confirms the use of the kurtosis to approximate the PDF of the sum of two independent GGRV by another GGD to get a good tail approximation, this may also confirms that the kurtosis measure the heavy tail. Unlike this observation, the $\gamma$ obtained by minimizing the CDF error is a little bit far from outcomes of the kurtosis method. To conclude, these three methods can be used according to the situation we are facing.

\subsection{PDF and CDF Simulations}
\begin{figure}[h]
\begin{center}
\includegraphics[width=0.95\columnwidth]{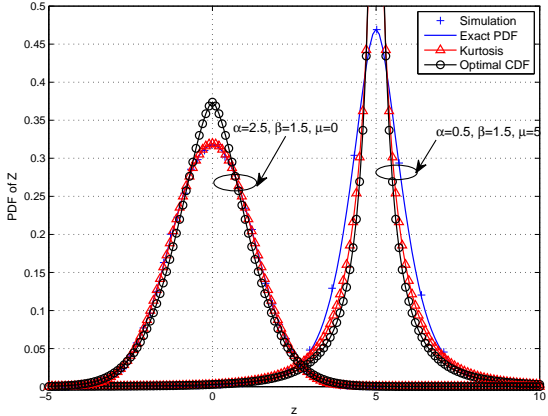}
\caption{Exact and approximated PDF of the sum of two GGRV, for $\beta=1.5$, $\delta=2$, $\sigma_1=1$, and two values of $\alpha$.}
\label{figure2}
\end{center}
\end{figure}
The illustrations in this section are made for $\beta=1.5$, $\delta=2$, and $\sigma_1=1$. In \figref{figure2}, the PDF of the sum distribution is drawn for two values of $\alpha$ ($0.5$ and $2.5$) and $\mu$ takes two values to split the curves of both cases. The exact and simulated PDF of $Z$ are drawn in the same figure among with the approximated PDF. The latter is computed using the kurtosis and the optimal CDF methods. It is clear that the exact PDF matchs perfectly the simulated PDF. Far from the mean, the approximated PDF appears close to the exact PDF and both methods have a good tail approximation. For $\alpha<2$, the kurtosis and the optimal CDF method are close to each other and match only the exact PDF at the tail with huge difference at the mean as mentioned in \cite{qian_sumgg}. However, for $\alpha>2$, the kurtosis method presents a good approximation of the PDF even around the mean, while the optimal CDF method represents a good approximation of the CDF as it will be seen later. We omit the optimal tail method here because it is close to the kurtosis method as shown in Table \ref{table1}. However one can draw it easily using the values available in Table \ref{table1}.

\begin{figure}[h]
\begin{center}
\includegraphics[width=0.95\columnwidth]{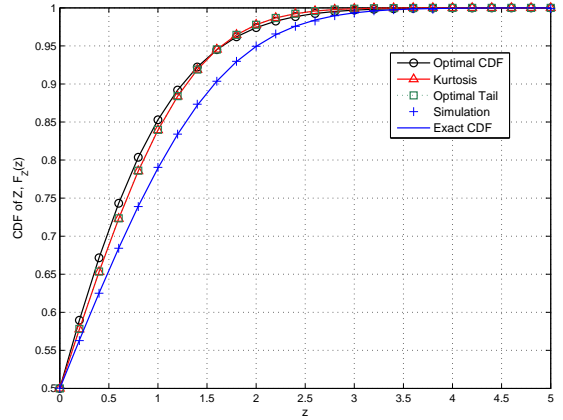}
\caption{CDF of the sum of two GGRV using the Kurtosis, optimal tail, and optimal CDF approximations methods, for $\alpha=2.5$, $\beta=1.5$, $\delta=2$, and $\sigma_1=1$}
\label{figure3}
\end{center}
\end{figure}

Our second illustration, highlighted in \figref{figure3}, consists of drawing the CDF of the sum for $\alpha=2.5$, $\beta=1.5$, $\delta=2$, and $\sigma_1=1$ using all three methods to approximate the PDF in linear scale. It is noticed that the results obtained from the optimal tail method (for $n=3$) are very close to those issued from the kurtosis method. Another observation is that all the methods are close to each other and close to the CDF at the saturation region, i.e. $F_Z(z)\approx 1$. This result is more detailed in the next figure which shows the complementary CDF in Log scale.

\begin{figure}[h]
\begin{center}
\includegraphics[width=0.95\columnwidth]{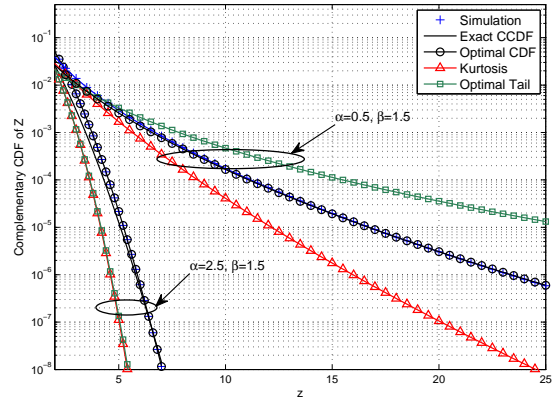}
\caption{Complementary CDF of the sum of two GGRV for two values of $\alpha$.}
\label{figure4}
\end{center}
\end{figure}
In \figref{figure4}, the complementary CDF of the distribution of the sum is drawn for two values of $\alpha$ ($2.5$ and $0.5$). For both cases, the approximated CCDF using the optimal CDF method matches the exact CCDF. However, as seen in \figref{figure2}, for $\alpha<2$, the CCDF obtained from the kurtosis and optimal tail methods is not too close to the exact CCDF. While, for $\alpha>1$, they are close to each other and asymptotically close to the exact CCDF.

\section{Conclusion}
The statistics of the distribution of the sum of two independent GGRV were derived in closed form in terms of the FHF and BFHF respectively, and the distribution of sum was approximated by another GGD using three estimation methods to get the shape parameter, namely kurtosis, optimal tail and optimal CDF depends on the treated case. Finally, all the methods present a good estimation of the shape parameter especially at the tail region.

\bibliography{IEEEfull,sum_generalizedG}
\bibliographystyle{IEEEtran}

\end{document}